\documentclass[letterpaper, 10pt, conference]{ieeeconf}      

\IEEEoverridecommandlockouts                              
\usepackage{graphics} 
\usepackage{comment}
\usepackage{epsfig} 
\usepackage{mathptmx} 
\usepackage{times} 
\usepackage{amsmath} 
\usepackage{amssymb}  
\usepackage{mathtools}
\usepackage{xcolor}
\usepackage{subfigure}
\usepackage{url}

\renewcommand{\baselinestretch}{0.91} 
\newtheorem{theorem}{Theorem}

\newtheorem{assumption}{Assumption}

\title{\LARGE \bf
Sliding-Mode Nash Equilibrium Seeking for a Quadratic Duopoly Game}

\author{Victor Hugo Pereira Rodrigues$^{1}$, Tiago Roux Oliveira$^{1}$, Miroslav Krsti{\' c}$^{2}$, Tamer Ba{\c s}ar$^{3}$.
\thanks{$^{1}$V. H. P. Rodrigues and T. R. Oliveira are with the Dept. of Eletronics and Telecommunication Engineering (DETEL), State University of Rio de Janeiro (UERJ), Rio de Janeiro -- RJ, Brazil. Email: {\tt\small victor.rodrigues@eng.uerj.br, tiagoroux@uerj.br}}%
       \thanks{$^{2}$M. Krsti{\' c} is with the Dept. of Mechanical and Aerospace Engineering, University of California -- San Diego (UCSD), La Jolla -- CA, USA. Email: {\tt\small krstic@ucsd.edu}} 
         \thanks{$^{3}$T. Ba{\c s}ar is with the Dept. of Electrical and Computer Engineering, University of Illinois Urbana-Champaign, Urbana -- IL, USA. Email: {\tt\small basar1@illinois.edu}}%
}

\begin{document}

\maketitle
\thispagestyle{empty}
\pagestyle{empty}

\begin{abstract}
This paper introduces a new method to achieve stable convergence to Nash equilibrium in duopoly noncooperative games. Inspired by the recent fixed-time Nash Equilibrium seeking (NES) \cite{PKB:2023} as well as prescribed-time extremum seeking (ES) \cite{YK:2022} and source seeking \cite{TK:2023} schemes, our approach employs a distributed sliding mode control (SMC) scheme, integrating extremum seeking with sinusoidal perturbation signals to estimate the pseudogradients of quadratic payoff functions. Notably, this is the first attempt to address noncooperative games without relying on models, combining classical extremum seeking with relay components instead of proportional control laws. We prove finite-time convergence of the closed-loop average system to Nash equilibrium using stability analysis techniques such as time-scaling, Lyapunov's direct method, and averaging theory for discontinuous systems. Additionally, we quantify the size of residual sets around the Nash equilibrium and validate our theoretical results through simulations.
%
\end{abstract}

\section{\textcolor{black}{INTRODUCTION}}

Game theory provides a conceptual framework for examining social interactions among competitive participants, utilizing mathematical models to grasp strategic decision-making \cite{Tirole:1991,BZ:2018}. Its utility extends across multiple domains such as engineering systems, biological behaviors, and financial markets, rendering it indispensable across diverse disciplines \cite{Basar:2019,Sastry:2013,BZ2:2018}. Considerable research has been conducted on differential games, encompassing both theoretical exploration and practical applications. 


Games are commonly categorized as cooperative and noncooperative \cite{Basar:1999}. In cooperative games, players establish enforceable agreements, whereas noncooperative games center on individual player actions and Nash equilibria \cite{N:1951}. Nash equilibrium, a pivotal concept in noncooperative game theory, characterizes a situation where no player can improve its payoff independently of other players in the game \cite{Basar:1999}.


Efforts towards achieving convergence to Nash equilibrium have persisted for decades \cite{LB:1987,B:1987}, with research delving into learning-based update schemes \cite{ZTB:2013}. Recent investigations have delved into extremum seeking (ES) methods for real-time computation of Nash equilibria in noncooperative games \cite{FKB:2012}, facilitating stable convergence in the absence of any model information \cite{KW:2000}.


From the literature of ES based on sliding modes, we can mention the publication from Özgüner and coauthors \cite{Ozguner_SMC}, where a periodic switching function was introduced to provide a state-feedback control approach. Following that, Oliveira et al. \cite{Oliveira_SMC1,Oliveira_SMC2} have expanded the application of this previous methodology by considering only output feedback. Notably, the same team has also introduced a different sliding mode-based ES controller using monitoring functions \cite{Aminde_SMC1,Aminde_SMC2}. However, none of such ES control strategies followed the original and simplest structure from the scheme with sinusoidal periodic perturbations \cite{KW:2000}.


On the other hand, recently some interesting fixed-time results for Nash Equilibrium seeking \cite{PKB:2023} (time-invariant, non-smooth) and prescribed-time schemes (time-varying, smooth) for ES \cite{YK:2022} and source seeking \cite{TK:2023} have been provided. In particular, the paper \cite{PKB:2023} proposes control laws which seem to be similar to those in higher-order SMC. For instance, when considering second-order SMC, the most known is the super-twisting algorithm (STA) \cite{MO:2012}. Indeed, there is an entire family of finite-time (or fixed-time) SMC schemes \cite{M:2011} which could be explored for extremum seeking. In this context, bringing the idea of a variable structure/sliding mode NES in finite (not fixed) time to the attention of the large SMC community is still missing and we believe has an important value. We have decided to start with the simplest (and more popular) case: first-order
SMC (FOSMC), simply a relay function. FOSMC is discontinuous but can assure finite-time convergence.

To bridge the existing gaps, our paper extends the classical ES algorithm \cite{KW:2000} by incorporating relay components in place of the usual proportional control laws with respect to the {\it pseudogradient}\footnote{The \textit{pseudogradient} is essentially a generalized notion of gradient that provides directional information allowing the analysis of equilibrium points in complex strategic interactions of concave $N$-players game. 
Hence, given a vector-valued function $J(\theta): \mathbb{R}^{N} \to \mathbb{R}^{N}$, we define its \textit{pseudogradient} $G$ as the vector-valued function $G(\theta):=\left[\frac{\partial J_1}{\partial \theta_1}\,, \frac{\partial J_2}{\partial \theta_2}\,, \ldots\,, \frac{\partial J_N}{\partial \theta_N}\right]^T$. For more details, see \cite[Eq. (3.9)]{R:1965}.} estimates in order to derive the first sliding mode version of the classical ES algorithm based on sinusoidal periodic perturbations.

We address the NES problem in two-player noncooperative games, employing the proposed distributed variable structure ES controller. Despite the inherent challenges of decentralized control in games, our approach ensures convergence to Nash equilibrium when both players adopt the proposed ES algorithms via sliding modes. We provide theoretical analysis for duopoly games, leveraging time-scaling techniques, constructing a Lyapunov function, and employing averaging methods to ensure closed-loop stability. Numerical simulations demonstrate the efficacy of our approach, highlighting its superiority over classical ES with linear proportional control laws \cite{KW:2000} since now the average closed-loop system is finite-time stable rather than only exponentially stable.



\section{Duopoly Game with Quadratic Payoffs: General Formulation}

In a duopoly game, the optimal outcomes for players P1 and P2, denoted by $y_{1}(t) \in \mathbb{R}$ and $y_{2}(t) \in \mathbb{R}$, respectively, are not solely determined by their individual actions or decision strategies, represented by $\theta_{1}(t)\in \mathbb{R}$ and $\theta_{2}(t)\in \mathbb{R}$, but also by their opponents'. Hence, the pair of decisions,
$\theta=[\theta_{1}\,,\theta_{2}]^T \in \mathbb{R}^{2}$, influences the payoff functions $J_{1}$ and $J_{2}$ of the two players. An optimal choice ($\theta(t)\equiv\theta^{*}$), called Nash equilibrium \cite{Tirole:1991}, captures the situation where no player benefits in terms of improvement of its own payoff by deviating unilaterally from its own policy. 

In particular, we examine duopoly games where each player's payoff function takes a quadratic form. This can be expressed for each player as a strictly concave function of its individual action:
\begin{align}
J_{1}(\theta(t)) &= \frac{H_{11}^{1}}{2}\theta_{1}^{2}(t) + \frac{H_{22}^{1}}{2}\theta_{2}^{2}(t) +H_{12}^{1}\theta_{1}(t)\theta_{2}(t) \nonumber \\
&\quad + h_{1}^{1}\theta_{1}(t) + h_{2}^{1}\theta_{2}(t) + c_{1}, \label{eq:J1} \\
J_{2}(\theta(t)) &= \frac{H_{11}^{2}}{2}\theta_{1}^{2}(t) + \frac{H_{22}^{2}}{2}\theta_{2}^{2}(t)+ H_{21}^{2}\theta_{1}(t)\theta_{2}(t) \nonumber \\
&\quad+h_{1}^{2}\theta_{1}(t) + h_{2}^{2}\theta_{2}(t) + c_{2}, \label{eq:J2}
\end{align}
Here, $J_{1}(\theta), J_{2}(\theta) : \mathbb{R}^{2} \to \mathbb{R}$ represent the payoff functions for players 1 and 2, respectively, while $\theta_{1}(t), \theta_{2}(t) \in \mathbb{R}$ denote the players' decision variables, and $H_{jk}^{i}$, $h_{j}^{i}$, $c_{i} \in \mathbb{R}$ are constants, with $H_{ii}^{i} < 0$, for all $i, j, k \in \{1, 2\}$. In the analysis of the duopoly game, the individual payoff functions exhibit strict concavity only in self-action, deviating from concavity in both actions. This implies the absence of a quadratic term in the opponent's action, leading to $H^i_{jj} = 0$ for all $j\neq i$, alongside the general non-zero nature of all other coefficients except for $c_i$, which are typically absent.

For the sake completeness, let us define the mathematical representation of the Nash equilibrium $\theta^{\ast}=[\theta^{\ast}_1\,, \theta^{\ast}_2]^T$ in a two-player game as below, for all $\theta_1\,, \theta_2$:
\begin{align} \label{Nashcu}
J_1(\theta_1^{\ast}\,,\theta_{2}^{\ast}) \geq J_1(\theta_1\,,\theta_{2}^{\ast}) ~~  \mbox{and} ~~~ J_2(\theta_1^{\ast}\,,\theta_{2}^{\ast}) \geq J_2(\theta_1^{\ast}\,,\theta_{2}). \!\!
\end{align}

Thus, no player has an incentive to independently deviate from $\theta^*$. In the duopoly context, $\theta_1$ and $\theta_2$ are elements of $\mathbb{R}$, representing the set of real numbers. 
In order to determine the Nash equilibrium solution in strictly concave quadratic games involving two players, where each action set is the entire real line, one would differentiate $J_{i}$ with respect to $\theta_{i}(t) \,, \forall i=1,2$, setting the resulting expressions equal to zero, and solving the set of equations thus obtained. This set of equations, which also provides a sufficient condition due to the strict concavity \cite{R:1965}, is
\begin{align}
\begin{cases}
H_{11}^{1}\theta_{1}^{*}+H_{12}^{1}\theta_{2}^{*}+h_{1}^{1}=0 \\ 
H_{21}^{2}\theta_{1}^{*}+H_{22}^{2}\theta_{2}^{*}+h_{2}^{2}=0 \\ 
\end{cases}\,,
\label{eq:NE_v0}
\end{align} 
which can be written in the matrix form as 
\begin{align}
\begin{bmatrix}
 H_{11}^{1} & H_{12}^{1} \\
 H_{21}^{2} & H_{22}^{2}
\end{bmatrix}
\begin{bmatrix}
\theta_{1}^{*} \\
\theta_{2}^{*} 
\end{bmatrix}
=-
\begin{bmatrix}
h_{1}^{1} \\
h_{2}^{2}  
\end{bmatrix}
\,.
\end{align}
Defining the Hessian matrix $H$, and vectors $\theta^*$ and $h$ by 
\begin{align}
H:=
\begin{bmatrix}
 H_{11}^{1} &  H_{12}^{1} \\
 H_{21}^{2} &  H_{22}^{2}
\end{bmatrix}
\,, \quad 
\theta^{*}:=
\begin{bmatrix}
\theta_{1}^{*} \\
\theta_{2}^{*} 
\end{bmatrix}
\,, \quad
h:=
\begin{bmatrix}
h_{1}^{1} \\
h_{2}^{2}  
\end{bmatrix}
\,, \label{eq:Htheta*h}
\end{align}
there exists a unique Nash Equilibrium at  $\theta^{*}\!=\!-H^{-1}h$, if $H$ is invertible (which holds in any realistic model of the duopoly game \cite[Chapter 4]{Basar:1999}):
\begin{align}
\begin{bmatrix}
\theta_{1}^{*} \\
\theta_{2}^{*} 
\end{bmatrix}
=-
\begin{bmatrix}
 H_{11}^{1} & H_{12}^{1} \\
 H_{21}^{2} &  H_{22}^{2}
\end{bmatrix}^{-1}
\begin{bmatrix}
h_{1}^{1} \\
h_{2}^{2}  
\end{bmatrix} \label{eq:NE_v2} \,.
\end{align}
%

\subsection{Nash Equilibrium Seeking}

The primary goal of this paper is to develop distributed variable structure policies via sliding modes to obtain NES in noncooperative duopoly games. Throughout our analysis,
it is assumed that $H_{ii}^{i}<0$ for both $i$ (a natural consequence of strict concavity), the matrix $H$ in (\ref{eq:Htheta*h}) has full rank, and all the parameters are unknown.

\begin{figure}[ht]
\begin{center}
\includegraphics[width=3.35 in]{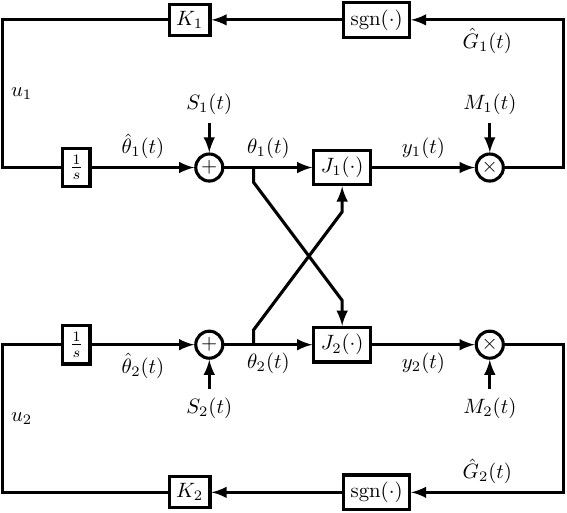}
\end{center}
\caption{Block diagram illustrating the NES strategy through distributed event-triggered control policies used by each player.}
\label{fig:blockDiagram_v3}
\end{figure}
The schematic diagram depicted in Fig.~\ref{fig:blockDiagram_v3} provides an overview of the proposed NES policy for each player, detailing their respective outputs as follows:
\begin{align}
\begin{cases}
y_{1}(t)&=J_{1}(\theta(t))\,, \\ 
y_{2}(t)&=J_{2}(\theta(t))\,.
\end{cases}
\label{eq:yi_v0}
\end{align}
The probing signals are
\begin{align}
\begin{cases}
S_{1}(t)&=a_{1}\sin(\omega_{1}t) \\
S_{2}(t)&=a_{2}\sin(\omega_{2}t)
\end{cases}
\,,\label{eq:Si}
\end{align}
and the demodulating signals are represented by
\begin{align}
\begin{cases}
M_{1}(t)&=\frac{2}{a_{1}}\sin(\omega_{1}t) \\
M_{2}(t)&=\frac{2}{a_{2}}\sin(\omega_{2}t)
\end{cases}
\,, \label{eq:Mi}
\end{align}
having non-zero constant amplitudes $a_{1}$ and $a_{2}$, both greater than zero, and frequencies $\omega_1$ and $\omega_2$, satisfying $\omega_1 \neq \omega_2$. These probing frequencies $\omega_i$ can be chosen as
\begin{equation}\label{omegadefinition}
\omega_i=\omega_{i}'\omega=\mathcal{O}(\omega)\,, \quad
i=1 \mbox{ or }2\,,
\end{equation}
where $\omega$ represents a positive constant and $\omega_{i}'$ is a rational number. One potential selection is detailed in \cite{GKN:2012}.

If we denote by $\hat{\theta}_{1}(t)$ and $\hat{\theta}_{2}(t)$ the estimates of $\theta^{\ast}_{1}$ and $\theta^{\ast}_{2}$, respectively, we can define the ``estimation error'' as follows:
\begin{align}
\begin{cases}
\tilde{\theta}_{1}(t)=\hat{\theta}_{1}(t)-\theta_{1}^{*}\\
\tilde{\theta}_{2}(t)=\hat{\theta}_{2}(t)-\theta_{2}^{*}
\end{cases}\,.\label{eq:tildeThetai}
\end{align}
The estimation of the {\it pseudogradient} of the unknown payoff functions 
\begin{align}
\begin{cases}
\hat{G}_{1}(t)&=M_1(t)y_1(t) \\
\hat{G}_{2}(t)&=M_2(t)y_2(t)
\end{cases}
\,, \label{eq:Gi_gold}
\end{align}
can be rewritten as
%
\begin{align}
\hat{G}_{1}(t)&= H^{1}_{11}\theta^{\ast}_{1} +H^{1}_{12}\theta^{\ast}_{2} + h^{1}_{1} \nonumber \\
&\quad +\mathcal{H}_{11}(t)\tilde{\theta}_{1}(t) +\mathcal{H}_{12}(t)\tilde{\theta}_{2}(t) \nonumber \\
&\quad+\frac{H^{1}_{11}}{a_{1}}\sin(\omega_{1}t)\tilde{\theta}_{1}^2(t)+ \frac{H^{1}_{22}}{a_{1}}\sin(\omega_{1}t)\tilde{\theta}_{2}^2(t)  \nonumber \\
&\quad+ \frac{2H^{1}_{12}}{a_{1}}\sin(\omega_{1}t)\tilde{\theta}_{1}(t)\tilde{\theta}_{2}(t) +\delta_{1}(t)\,, \label{eq:hatG1_20240318} \\
\hat{G}_{2}(t)&= H^{2}_{21}\theta^{\ast}_{1} +H^{2}_{22}\theta^{\ast}_{2} + h^{2}_{2} \nonumber \\
&\quad +\mathcal{H}_{21}(t)\tilde{\theta}_{1}(t) +\mathcal{H}_{22}(t)\tilde{\theta}_{2}(t) \nonumber \\
&\quad+\frac{H^{2}_{11}}{a_{2}}\sin(\omega_{2}t)\tilde{\theta}_{1}^2(t)+ \frac{H^{2}_{22}}{a_{2}}\sin(\omega_{2}t)\tilde{\theta}_{2}^2(t)  \nonumber \\
&\quad+ \frac{2H^{2}_{21}}{a_{2}}\sin(\omega_{2}t)\tilde{\theta}_{1}(t)\tilde{\theta}_{2}(t) +\delta_{2}(t)\,, \label{eq:hatG2_20240318}
\end{align}
with the time-varying parameters given by
\begin{align}
\mathcal{H}_{11}(t)&=H^{1}_{11}-H^{1}_{11}\cos(2\omega_{1}t)   \nonumber \\
&\quad+\frac{2h^{1}_{1}}{a_{1}}\sin(\omega_{1}t)+ \frac{2H^{1}_{11}\theta^{\ast}_{1}}{a_{1}}\sin(\omega_{1}t) \nonumber \\
&\quad +\frac{H^{1}_{12}a_{2}}{a_{1}}\cos[(\omega_{1}+\omega_{2})t] \nonumber \\
&\quad +\frac{2H^{1}_{12}\theta^{\ast}_{2}}{a_{1}}\sin(\omega_{1}t)+\frac{H^{1}_{12}a_{2}}{a_{1}}\cos[(\omega_{1}-\omega_{2})t] \,, \\
\mathcal{H}_{12}(t)&= H^{1}_{12}- H^{1}_{12}\cos(2\omega_{1}t) \nonumber \\
&\quad+ \frac{2h^{1}_{2}}{a_{1}}\sin(\omega_{1}t)+ \frac{2H^{1}_{12}\theta^{\ast}_{1}}{a_{1}}\sin(\omega_{1}t) \nonumber \\
&\quad + \frac{2H^{1}_{22}\theta^{\ast}_{2}}{a_{1}}\sin(\omega_{1}t)+\frac{H^{1}_{22}a_{2}}{a_{1}}\cos[(\omega_{1}-\omega_{2})t] \nonumber \\
&\quad+\frac{H^{1}_{22}a_{2}}{a_{1}}\cos[(\omega_{1}+\omega_{2})t]\,, \\
\mathcal{H}_{21}(t)&= H^{2}_{21}- H^{2}_{21}\cos(2\omega_{2}t) \nonumber \\
&\quad+ \frac{2h^{2}_{1}}{a_{2}}\sin(\omega_{2}t)+ \frac{2H^{2}_{21}\theta^{\ast}_{1}}{a_{2}}\sin(\omega_{2}t) \nonumber \\
&\quad + \frac{2H^{2}_{22}\theta^{\ast}_{2}}{a_{2}}\sin(\omega_{2}t)+\frac{H^{2}_{22}a_{1}}{a_{2}}\cos[(\omega_{2}-\omega_{1})t] \nonumber \\
&\quad+\frac{H^{2}_{22}a_{1}}{a_{2}}\cos[(\omega_{2}+\omega_{1})t]\,, \\
\mathcal{H}_{22}(t)&=H^{2}_{22}-H^{2}_{22}\cos(2\omega_{2}t)   \nonumber \\
&\quad+\frac{2h^{2}_{2}}{a_{2}}\sin(\omega_{2}t)+ \frac{2H^{2}_{22}\theta^{\ast}_{2}}{a_{2}}\sin(\omega_{2}t) \nonumber \\
&\quad +\frac{H^{2}_{21}a_{1}}{a_{2}}\cos[(\omega_{2}+\omega_{1})t] \nonumber \\
&\quad +\frac{2H^{2}_{21}\theta^{\ast}_{1}}{a_{2}}\sin(\omega_{2}t)+\frac{H^{2}_{21}a_{1}}{a_{2}}\cos[(\omega_{2}-\omega_{1})t] \,,
\end{align}
\vspace{-1cm}
\begin{small}
\begin{align}
\delta_{1}(t)&= \frac{H^{1}_{11}a_{1}}{2}\sin(\omega_{1}t)+\frac{H^{1}_{11}a_{1}}{4}\cos(\omega_{1}t)+\frac{H^{1}_{11}a_{1}}{4}\cos(3\omega_{1}t) \nonumber \\ 
&\quad +\frac{H^{1}_{22}a_{2}^2}{2a_{1}}\sin(\omega_{1}t)-\frac{H^{1}_{22}a_{2}^2}{2a_{1}}\cos[(\omega_{1}-2\omega_{2})t] \nonumber \\
&\quad+\frac{H^{1}_{22}a_{2}^2}{2a_{1}}\cos[(\omega_{1}+2\omega_{2})t]+\frac{2H^{1}_{12}\theta^{\ast}_{1}\theta^{\ast}_{2}}{a_{1}}\sin(\omega_{1}t)   \nonumber \\
&\quad-h^{1}_{1}\cos(2\omega_{1}t)  \nonumber \\
&\quad+ \frac{a_{2}h^{1}_{2}}{a_{1}}\cos[(\omega_{1}-\omega_{2})t]+ \frac{a_{2}h^{1}_{2}}{a_{1}}\cos[(\omega_{1}+\omega_{2})t] \nonumber \\
&\quad-H^{1}_{11}\theta^{\ast}_{1}\cos(2\omega_{1}t)  \nonumber \\
&\quad + \frac{H^{1}_{12}a_{2}\theta^{\ast}_{1}}{a_{1}}\cos[(\omega_{1}-\omega_{2})t]+ \frac{H^{1}_{12}a_{2}\theta^{\ast}_{1}}{a_{1}}\cos[(\omega_{1}+\omega_{2})t]\nonumber \\
&\quad -H^{1}_{12}\theta^{\ast}_{2}\cos(2\omega_{1}t)  \nonumber \\
&\quad + \frac{H^{1}_{22}a_{2}\theta^{\ast}_{2}}{a_{1}}\cos[(\omega_{1}-\omega_{2})t] - \frac{H^{1}_{22}a_{2}\theta^{\ast}_{2}}{a_{1}}\cos[(\omega_{1}+\omega_{2})t]\nonumber \\
&\quad+\frac{2c_{1}}{a_{1}}\sin(\omega_{1}t) + \frac{2h^{1}_{1}\theta^{\ast}_{1}}{a_{1}}\sin(\omega_{1}t) + \frac{2h^{1}_{2}\theta^{\ast}_{2}}{a_{1}}\sin(\omega_{1}t)  \nonumber \\
&\quad+\frac{H^{1}_{11}\theta^{\ast 2}_{1}}{a_{1}}\sin(\omega_{1}t) + \frac{H^{1}_{22}\theta^{\ast 2}_{2}}{a_{1}}\sin(\omega_{1}t) \nonumber \\
&\quad +\frac{H^{1}_{12}a_{2}}{2}\cos(\omega_{2}t)-\frac{H^{1}_{12}a_{2}}{2}\cos[(2\omega_{1}-\omega_{2})t] \nonumber \\
&\quad+\frac{H^{1}_{12}a_{2}}{2}\cos(\omega_{2}t)+\frac{H^{1}_{12}a_{2}}{2}\cos[(2\omega_{1}+\omega_{2})t]\,,
\end{align}
\end{small}

\begin{small}
\begin{align}
\delta_{2}(t)&= \frac{H^{2}_{11}a_{2}}{2}\sin(\omega_{2}t)+\frac{H^{2}_{11}a_{2}}{4}\cos(\omega_{2}t)+\frac{H^{2}_{11}a_{2}}{4}\cos(3\omega_{2}t) \nonumber \\ 
&\quad +\frac{H^{2}_{22}a_{2}}{2}\sin(\omega_{2}t)-\frac{H^{2}_{22}a_{2}}{2}\cos[(\omega_{2}-2\omega_{1})t] \nonumber \\
&\quad+\frac{H^{2}_{22}a_{2}}{2}\cos[(\omega_{2}+2\omega_{1})t]+\frac{2H^{2}_{21}\theta^{\ast}_{1}\theta^{\ast}_{2}}{a_{2}}\sin(\omega_{2}t)   \nonumber \\
&\quad-h^{2}_{2}\cos(2\omega_{2}t)  \nonumber \\
&\quad+ h^{1}_{2}\cos[(\omega_{2}-\omega_{1})t]+ \frac{h^{2}_{2}}{a_{2}}\cos[(\omega_{1}+\omega_{2})t] \nonumber \\
&\quad-H^{2}_{11}\theta^{\ast}_{1}\cos(2\omega_{2}t)  \nonumber \\
&\quad + \frac{H^{2}_{21}a_{2}\theta^{\ast}_{1}}{a_{1}}\cos[(\omega_{2}-\omega_{1})t]+ H^{2}_{21}\theta^{\ast}_{1}\cos[(\omega_{1}+\omega_{2})t]\nonumber \\
&\quad -H^{2}_{21}\theta^{\ast}_{2}\cos(2\omega_{2}t)  \nonumber \\
&\quad + \frac{H^{2}_{22}a_{1}\theta^{\ast}_{2}}{a_{2}}\cos[(\omega_{2}-\omega_{1})t] - \frac{H^{2}_{22}a_{1}\theta^{\ast}_{2}}{a_{2}}\cos[(\omega_{2}+\omega_{1})t]\nonumber \\
&\quad+\frac{2c_{2}}{a_{2}}\sin(\omega_{2}t) + \frac{2h^{2}_{1}\theta^{\ast}_{1}}{a_{2}}\sin(\omega_{2}t) + \frac{2h^{2}_{2}\theta^{\ast}_{2}}{a_{2}}\sin(\omega_{2}t)  \nonumber \\
&\quad+\frac{H^{2}_{11}\theta^{\ast 2}_{1}}{a_{2}}\sin(\omega_{2}t) + \frac{H^{2}_{22}\theta^{\ast 2}_{2}}{a_{2}}\sin(\omega_{2}t) \nonumber \\
&\quad +\frac{H^{2}_{21}a_{1}}{2}\cos(\omega_{1}t)-\frac{H^{2}_{21}a_{1}}{2}\cos[(2\omega_{2}-\omega_{1})t] \nonumber \\
&\quad+\frac{H^{2}_{21}a_{1}}{2}\cos(\omega_{1}t)+\frac{H^{2}_{21}a_{1}}{2}\cos[(2\omega_{2}+\omega_{1})t]\,.
\end{align}
\end{small}
It is worth noting that the quadratic terms in (\ref{eq:hatG1_20240318}) and (\ref{eq:hatG2_20240318})  can be disregarded in a local analysis \cite{K:2014}. Moreover, from (\ref{eq:NE_v0}), $H_{11}^{1}\theta_{1}^{*}+H_{12}^{1}\theta_{2}^{*}+h_{1}^{1}=0$ and $H_{21}^{2}\theta_{1}^{*}+H_{22}^{2}\theta_{2}^{*}+h_{2}^{2}=0$. Therefore, the {\it pseudogradient} estimates are locally expressed as
\begin{align}
\hat{G}_{1}(t)&= \mathcal{H}_{11}(t)\tilde{\theta}_{1}(t) +\mathcal{H}_{12}(t)\tilde{\theta}_{2}(t) +\delta_{1}(t)\,, \label{eq:hatG1_v2_20240318} \\
\hat{G}_{2}(t)&= \mathcal{H}_{21}(t)\tilde{\theta}_{1}(t) +\mathcal{H}_{22}(t)\tilde{\theta}_{2}(t) +\delta_{2}(t)\,. \label{eq:hatG2_v2_20240318}
\end{align}

Now, considering the time-derivative of (\ref{eq:tildeThetai}) and referring to the NES scheme illustrated in Fig.~\ref{fig:blockDiagram_v3}, we can derive the dynamics governing $\hat{\theta}(t)$ and $\tilde{\theta}(t)$ as follows:
\begin{align}
\frac{d\tilde{\theta}_{1}(t)}{dt}&=\frac{d\hat{\theta}_{1}(t)}{dt}=u_{1}(t) \label{eq:dotTildeTheta_1}\,, \\
\frac{d\tilde{\theta}_{2}(t)}{dt}&=\frac{d\hat{\theta}_{2}(t)}{dt}=u_{2}(t)\,, \label{eq:dotTildeTheta_2}
\end{align}
where $u_{1}(t)$ and $u_{2}(t)$ represent the NES control laws to be formulated.

\subsection{Distributed Variable Structure Actions for Nash Equilibrium Seeking via Sliding Modes}  

In the linear approach, the continuous-time feedback law structures remain fixed as 
\begin{align}
\begin{cases}
u_{1}(t)&=K_{1}\hat{G}_{1}(t) \\
u_{2}(t)&=K_{2}\hat{G}_{2}(t) 
\end{cases}\,,
\quad \text{for all } t\geq 0\,. \label{eq:U_continuous}
\end{align}
with the constants $K_{1}$ and $K_{2}$ chosen such that 
\begin{align}
K=\begin{bmatrix} K_{1} & 0 \\
                    0   & K_{2} \end{bmatrix}\,,
\end{align}
ensuring that $KH$ is Hurwitz.

In variable structure systems (VSS), the control is allowed to dynamically change its structure, switching between various continuous functions of the state. The challenge in variable structure design lies in selecting the parameters for each structure and defining the switching logic. This concept of structural flexibility enables the creation of novel system behaviors by combining segments of trajectories from various structures. Additionally, a fundamental aspect of variable structure systems is the ability to generate trajectories that are not naturally present in any individual structure. These trajectories define a unique form of motion known as the sliding mode. 

In this paper, we design variable structure actions for Nash Equilibrium Seeking through the following sliding mode control laws:
\begin{align}
\begin{cases}
u_{1}&=K_{1}\text{sgn}(\hat{G}_{1}(t)) \\
u_{2}&=K_{2}\text{sgn}(\hat{G}_{2}(t))
\end{cases}\,. \label{eq:U_discontinuous}
\end{align}

\subsection{Closed-loop System}

Plugging the proposed sliding mode based control laws in (\ref{eq:U_discontinuous}) into the time-derivatives of (\ref{eq:hatG1_v2_20240318}), (\ref{eq:hatG2_v2_20240318}), (\ref{eq:dotTildeTheta_1}) and equation (\ref{eq:dotTildeTheta_2}) yields, for all $t\geq 0$, the following dynamics governing $\hat{G}(t)$ and $\tilde{\theta}(t)$:
\begin{align}
\frac{d\hat{G}_{1}(t)}{dt}
&= \mathcal{H}_{11}(t)K_{1}\text{sgn}(\hat{G}_{1}(t)) +\mathcal{H}_{12}(t)K_{2}\text{sgn}(\hat{G}_{2}(t))+\nonumber \\
&\quad+\frac{d\mathcal{H}_{11}(t)}{dt}\tilde{\theta}_{1}(t) +\frac{d\mathcal{H}_{12}(t)}{dt}\tilde{\theta}_{2}(t)+\frac{d\delta_{1}(t)}{dt}\,, \label{eq:dotHatG1_v1_20240319} \\
\frac{d\hat{G}_{2}(t)}{dt}
&= \mathcal{H}_{21}(t)K_{1}\text{sgn}(\hat{G}_{1}(t)) +\mathcal{H}_{22}(t)K_{2}\text{sgn}(\hat{G}_{2}(t))+\nonumber \\
&\quad+\frac{d\mathcal{H}_{21}(t)}{dt}\tilde{\theta}_{1}(t) +\frac{d\mathcal{H}_{22}(t)}{dt}\tilde{\theta}_{2}(t)+\frac{d\delta_{2}(t)}{dt}\,, \label{eq:dotHatG2_v1_20240319} \\
\frac{d\tilde{\theta}_{1}(t)}{dt}
&= K_{1}\text{sgn}(\hat{G}_{1}(t))\,, \label{eq:dotTildeTheta1_v1_20240319} \\
\frac{d\tilde{\theta}_{2}(t)}{dt}
&= K_{2}\text{sgn}(\hat{G}_{2}(t))\,. \label{eq:dotTildeTheta2_v1_20240319}
\end{align}
Hence, in a concise vector-form representation, we have
\begin{align}
\frac{d\hat{G}(t)}{dt}&= \mathcal{H}(t)K\text{sgn}(\hat{G}(t))+\frac{d\mathcal{H}(t)}{dt}\tilde{\theta}(t)+\frac{d\delta(t)}{dt}\,, \label{eq:dotHatG_v1_20240319} \\
\frac{d\tilde{\theta}(t)}{dt}&= K\text{sgn}(\hat{G}(t))\,, \label{eq:dotTildeTheta_v1_20240319} 
\end{align}
where 
\begin{align}
\text{sgn}(\hat{G}(t))&:=
							 \begin{bmatrix} 
									\text{sgn}(\hat{G}_{1}(t)) \\
									\text{sgn}(\hat{G}_{2}(t))
							 \end{bmatrix} \,, \quad
K:=\begin{bmatrix} K_{1} & 0 \\
                    0   & K_{2} \end{bmatrix}\,, \nonumber \\
\mathcal{H}(t)&:=\begin{bmatrix} 
									\mathcal{H}_{11}(t) & \mathcal{H}_{12}(t) \\
									\mathcal{H}_{21}(t) & \mathcal{H}_{22}(t)
							 \end{bmatrix} \,,
\mbox{ and } 
\delta(t):=\begin{bmatrix} 
									\delta_{1}(t) \\
									\delta_{2}(t)
							 \end{bmatrix} \,.
\end{align}
The closed-loop system described by (\ref{eq:dotHatG_v1_20240319}) and (\ref{eq:dotTildeTheta_v1_20240319}) highlights a crucial point: even if the product $\mathcal{H}(t)K$ on averaging sense would form a Hurwitz matrix, the convergence to the equilibrium $\hat{G}\equiv0$ and $\tilde{\theta}\equiv0$ would not be directly guaranteed through the
$\tilde{\theta}$--dynamics. 
On the other hand, it is important to note that the disturbances $\delta(t)$ and $\frac{d\delta(t)}{dt}$ as well as the time-varying matrix $\frac{d\mathcal{H}(t)}{dt}$ possess zero mean values, which will be helpful in the following averaging analysis.

\subsection{Time-Scaling Procedure} \label{sec:time_scaling}

To facilitate the stability analysis of the closed-loop system, we introduce a convenient time scale. By examining (\ref{omegadefinition}), it is evident that the dither frequencies (\ref{eq:Si}) and (\ref{eq:Mi}), as well as their combinations, are rational. Additionally, there exists a time period $T$ such that
\begin{align}
T&= 2\pi \times \text{LCM}\left\{\frac{1}{\omega_{1}},\frac{1}{\omega_{2}}\right\}\,,\label{eq:period_T}
\end{align}
where LCM denotes the least common multiple. Hence, we define a new time scale for the dynamics (\ref{eq:dotHatG_v1_20240319}) and (\ref{eq:dotTildeTheta_v1_20240319}) using the transformation $\bar{t}=\omega t$, where
\begin{align}
\omega&:=\frac{2\pi}{T}\,. \label{eq:new_omega}
\end{align}
Consequently, the system (\ref{eq:dotHatG_v1_20240319}) and (\ref{eq:dotTildeTheta_v1_20240319}) can be expressed as
\begin{align}
\frac{d\hat{G}(\bar{t})}{d\bar{t}}&= \frac{1}{\omega}\mathcal{H}(\bar{t})K\text{sgn}(\hat{G}(\bar{t}))+\frac{1}{\omega}\frac{d\mathcal{H}(\bar{t})}{d\bar{t}}\tilde{\theta}(\bar{t})+\frac{1}{\omega}\frac{d\delta(\bar{t})}{d\bar{t}}\,, \label{eq:dotHatG_v2_20240319} \\
\frac{d\tilde{\theta}(\bar{t})}{d\bar{t}}&= \frac{1}{\omega}K\text{sgn}(\hat{G}(\bar{t}))\,. \label{eq:dotTildeTheta_v2_20240319} 
\end{align}
Now, we can implement a suitable averaging mechanism within the transformed time scale $\bar{t}$ based on the dynamics (\ref{eq:dotHatG_v2_20240319}) and (\ref{eq:dotTildeTheta_v2_20240319}). Despite the discontinuity on the right-hand side of equations (\ref{eq:dotHatG_v2_20240319}) and (\ref{eq:dotTildeTheta_v2_20240319}), the closed-loop system maintains its periodicity over time due to the periodic probing and demodulation signals. This unique characteristic allows for the application of the averaging results established by Plotnikov \cite{P:1979} to this particular setup. 

\subsection{Average Closed-Loop System}

Let us define the augmented state as:
\begin{align}
X^{T}(\bar{t}):=\begin{bmatrix} \hat{G}(\bar{t})\,, \tilde{\theta}(\bar{t}) \end{bmatrix}\,,
\end{align}
which reduces the system (\ref{eq:dotHatG_v2_20240319}) and (\ref{eq:dotTildeTheta_v2_20240319}) to:
\begin{align}
\dfrac{dX(\bar{t})}{d\bar{t}}&=\dfrac{1}{\omega}\mathcal{F}\left(\bar{t},X,\dfrac{1}{\omega}\right)\,. \label{eq:dotX_event}
\end{align}

The system (\ref{eq:dotX_event}) features a small parameter $1/\omega$ and a $T$-periodic function $\mathcal{F}\left(\bar{t},X,\dfrac{1}{\omega}\right)$ in $\bar{t}$. Therefore, the averaging theorem \cite{P:1979} can be applied to $\mathcal{F}\left(\bar{t},X,\dfrac{1}{\omega}\right)$ at $\displaystyle \lim_{\omega\to \infty}\dfrac{1}{\omega}=0$, following the principles established by Plotnikov \cite{P:1979}.

By employing the averaging of (\ref{eq:dotX_event}), we derive the following average system:
\begin{align}
\dfrac{dX^{\rm{av}}(\bar{t})}{d\bar{t}}&=\dfrac{1}{\omega}\mathcal{F}^{\rm{av}}\left(X^{\rm{av}}\right) \,, \label{eq:dotXav_event_1} \\
\mathcal{F}^{\rm{av}}\left(X^{\rm{av}}\right)&=\dfrac{1}{T}\int_{0}^{T}\mathcal{F}\left(\xi,X^{\rm{av}},0\right)d\xi \label{eq:mathcalFav_event} \,,
\end{align}
where we ``freeze'' the average states of $\hat{G}(\bar{t})$ and $\tilde{\theta}(\bar{t})$, resulting in the following equations, for all $t \geq 0$:
\begin{align}
\frac{d\hat{G}^{\rm{av}}(\bar{t})}{d\bar{t}}&=\frac{1}{\omega}HK\text{sgn}(\hat{G}^{\rm{av}}(\bar{t}))\,, \label{eq:dotHatGav_event_1} \\
\frac{d\tilde{\theta}^{\rm{av}}(\bar{t})}{d\bar{t}}&=\frac{1}{\omega}K\text{sgn}(\hat{G}^{\rm{av}}(\bar{t}))\,, \label{eq:dotTildeThetaAv_event_1} \\
\hat{G}^{\rm{av}}(\bar{t})&= H\tilde{\theta}^{\rm{av}}(\bar{t})\,. \label{eq:hatGav_event_1}
\end{align}
Now, we introduce a necessary and sufficient condition regarding the matrix $HK$. 

\smallskip
\begin{assumption} \label{einsteinnewton}
The matrix $HK = \begin{bmatrix} H_{11}^{1}K_{1} & H_{12}^{1}K_{2} \\ H_{21}^{2}K_{1} & H_{22}^{2}K_{2}\end{bmatrix}$ is strictly diagonally dominant, {\it i.e.}, for $j \neq i$, 
\begin{align} \label{xuxaangelica}
|H_{ij}^{i}K_{j}| < |H_{ii}^{i}K_{i}|\,, \quad i \in \{1,2\}\,.
\end{align}
Additionally, 
\begin{align}
\text{sgn}(K_{i})=-\text{sgn}(H_{ii}^{i})\,. \quad i \in \{1,2\}\,. \label{eq:sgnKi}
\end{align}
\end{assumption} 
\smallskip

According to Assumption~\ref{einsteinnewton}, we can establish the existence and uniqueness of the Nash equilibrium. This conclusion stems from the fact that strictly diagonally dominant matrices, as proven by the Levy–Desplanques theorem \cite{taussky1949,horn1985}, are nonsingular. 

\section{Stability Analysis}

The next theorem guarantees finite-time stability of the proposed NES control system via sliding modes, as depicted in Fig.~\ref{fig:blockDiagram_v3}.

\begin{theorem} \label{thm:NETESC_2}
Consider the closed-loop average dynamics of the {\it pseudogradient} estimate (\ref{eq:dotHatGav_event_1}), and that Assumption~\ref{einsteinnewton} holds. For a sufficiently large $\omega>0$, defined in (\ref{eq:new_omega}), there exists a finite time $t_{s} \in \left(0\,,\dfrac{\|\hat{G}^{\text{av}}(0)\|}{\min\{|H_{11}^{1}K_{1}|-|H_{12}^{1}K_{2}|,|H_{22}^{2}K_{2}|-|H_{21}^{2}K_{1}|\}}\right]$ such that the sliding motion on $\hat{G}_{\rm{av}}(t)\equiv0$ occurs. Furthermore, for the non-average system (\ref{eq:dotTildeTheta_v1_20240319}), there exist constants $m$ and $M_{\theta}$ 
such that
\begin{align}
\|\theta(t)\mathbb{-}\theta^{\ast}\|\mathbb{\leq} \mathcal{O}\left(a\mathbb{+}\dfrac{1}{\omega}\right) + \begin{cases} M_{\theta}\|\theta(0) \mathbb{-} \theta^{\ast}\|\mathbb{-} m t, \forall t \in [0, t_{s}),  \\ 0, \forall t \in [t_{s},+\infty), \end{cases}\!\!\!\!\!\!\!\! \label{eq:normTheta_thm2} 
\end{align}
where $a$ is defined as $a=\sqrt{a_{1}^{2}+a_{2}^{2}}$, with $a_1$ and $a_2$ being specified in \eqref{eq:Si} and \eqref{eq:Mi}. The constant $m$ depends on both the control gain $K$ and the Hessian $H$, while $M_{\theta}$ solely relies on the Hessian $H$.
\end{theorem}

\textit{Proof:} Let us examine the following candidate Lyapunov function for the average system  (\ref{eq:dotHatGav_event_1}):
\begin{align}
V_{\text{av}}(\bar{t})&=[\hat{G}_{1}^{\text{av}}(\bar{t})]^2+[\hat{G}_{2}^{\text{av}}(\bar{t})]^2 \,. \label{eq:Vav_event_pf2}
\end{align} 
The time derivative of (\ref{eq:Vav_event_pf2}) is given by  
\begin{align}
&\frac{dV_{\text{av}}(\bar{t})}{d\bar{t}}=2\hat{G}_{1}^{\text{av}}(\bar{t})\frac{d\hat{G}_{1}^{\text{av}}(\bar{t})}{d\bar{t}}+2\hat{G}_{2}^{\text{av}}(\bar{t})\frac{d\hat{G}_{2}^{\text{av}}(\bar{t})}{d\bar{t}}\nonumber \\
&=+\frac{2}{\omega}\hat{G}_{1}^{\text{av}}(\bar{t})(H_{11}^{1}K_{1}\text{sgn}(\hat{G}_{1}^{\text{av}}(\bar{t})) +H_{12}^{1}K_{2}\text{sgn}(\hat{G}_{2}^{\text{av}}(\bar{t})))+ \nonumber \\
&\quad+\frac{2}{\omega}\hat{G}_{2}^{\text{av}}(\bar{t})(H_{21}^{2}K_{1}\text{sgn}(\hat{G}_{1}^{\text{av}}(\bar{t})) +H_{22}^{2}K_{2}\text{sgn}(\hat{G}_{2}^{\text{av}}(\bar{t}))) \nonumber \\
&=\frac{2}{\omega}H_{11}^{1}K_{1}\text{sgn}(\hat{G}_{1}^{\text{av}}(\bar{t}))\hat{G}_{1}^{\text{av}}(\bar{t}) +\frac{2}{\omega}H_{12}^{1}K_{2}\text{sgn}(\hat{G}_{2}^{\text{av}}(\bar{t}))\hat{G}_{1}^{\text{av}}(\bar{t})+ \nonumber \\
&\quad+\frac{2}{\omega}H_{21}^{2}K_{1}\text{sgn}(\hat{G}_{1}^{\text{av}}(\bar{t}))\hat{G}_{2}^{\text{av}}(\bar{t}) +\frac{2}{\omega}H_{22}^{2}K_{2}\text{sgn}(\hat{G}_{2}^{\text{av}}(\bar{t}))\hat{G}_{2}^{\text{av}}(\bar{t}) \nonumber \\
&=\frac{2}{\omega}H_{11}^{1}K_{1}|\hat{G}_{1}^{\text{av}}(\bar{t})| +\frac{2}{\omega}H_{12}^{1}K_{2}\text{sgn}(\hat{G}_{2}^{\text{av}}(\bar{t}))\hat{G}_{1}^{\text{av}}(\bar{t})+ \nonumber \\
&\quad+\frac{2}{\omega}H_{21}^{2}K_{1}\text{sgn}(\hat{G}_{1}^{\text{av}}(\bar{t}))\hat{G}_{2}^{\text{av}}(\bar{t}) +\frac{2}{\omega}H_{22}^{2}K_{2}|\hat{G}_{2}^{\text{av}}(\bar{t})|\,.
\label{eq:dotVav_event_1_pf2}
\end{align}
By using (\ref{eq:sgnKi}), equation (\ref{eq:dotVav_event_1_pf2}) becomes
\begin{align}
\frac{dV_{\text{av}}(\bar{t})}{d\bar{t}}&=-\frac{2}{\omega}|H_{11}^{1}K_{1}||\hat{G}_{1}^{\text{av}}(\bar{t})| +\frac{2}{\omega}H_{12}^{1}K_{2}\text{sgn}(\hat{G}_{2}^{\text{av}}(\bar{t}))\hat{G}_{1}^{\text{av}}(\bar{t})+ \nonumber \\
&\quad+\frac{2}{\omega}H_{21}^{2}K_{1}\text{sgn}(\hat{G}_{1}^{\text{av}}(\bar{t}))\hat{G}_{2}^{\text{av}}(\bar{t}) -\frac{2}{\omega}|H_{22}^{2}K_{2}||\hat{G}_{2}^{\text{av}}(\bar{t})|\,,
\label{eq:dotVav_event_1_pf2_20240514}
\end{align}
whose upper bound can be expressed as
\begin{small}
\begin{align}
\frac{dV_{\text{av}}(\bar{t})}{d\bar{t}}&\leq-\frac{2}{\omega}|H_{11}^{1}K_{1}||\hat{G}_{1}^{\text{av}}(\bar{t})| +\frac{2}{\omega}|H_{12}^{1}K_{2}||\text{sgn}(\hat{G}_{2}^{\text{av}}(\bar{t}))||\hat{G}_{1}^{\text{av}}(\bar{t})|+ \nonumber \\
&\quad+\frac{2}{\omega}|H_{21}^{2}K_{1}||\text{sgn}(\hat{G}_{1}^{\text{av}}(\bar{t}))||\hat{G}_{2}^{\text{av}}(\bar{t})| -\frac{2}{\omega}|H_{22}^{2}K_{2}||\hat{G}_{2}^{\text{av}}(\bar{t})| \nonumber \\
&=-\frac{2}{\omega}|H_{11}^{1}K_{1}||\hat{G}_{1}^{\text{av}}(\bar{t})| +\frac{2}{\omega}|H_{12}^{1}K_{2}||\text{sgn}(\hat{G}_{2}^{\text{av}}(\bar{t}))||\hat{G}_{1}^{\text{av}}(\bar{t})|+ \nonumber \\
&\quad+\frac{2}{\omega}|H_{21}^{2}K_{1}||\text{sgn}(\hat{G}_{1}^{\text{av}}(\bar{t}))||\hat{G}_{2}^{\text{av}}(\bar{t})| -\frac{2}{\omega}|H_{22}^{2}K_{2}||\hat{G}_{2}^{\text{av}}(\bar{t})| \nonumber \\
&=-\frac{2}{\omega}(|H_{11}^{1}K_{1}|-|H_{12}^{1}K_{2}|)|\hat{G}_{1}^{\text{av}}(\bar{t})|+ \nonumber \\
&\quad-\frac{2}{\omega}(|H_{22}^{2}K_{2}|-|H_{21}^{2}K_{1}|)|\hat{G}_{2}^{\text{av}}(\bar{t})|. \label{eq:dotVav_event_2_pf2}
\end{align}
\end{small}
$\!\!$From the strict diagonal dominance assumption, equation (\ref{xuxaangelica}), one can introduce positive constants $\alpha_{1}=|H_{11}^{1}K_{1}|-|H_{12}^{1}K_{2}|$ and $\alpha_{2}=|H_{22}^{2}K_{2}|-|H_{21}^{2}K_{1}|$ such that inequality (\ref{eq:dotVav_event_2_pf2}) is rewritten as
\begin{align}
\frac{dV_{\text{av}}(\bar{t})}{d\bar{t}}&\leq-\frac{2}{\omega}\alpha_{1}|\hat{G}_{1}^{\text{av}}(\bar{t})|-\frac{2}{\omega}\alpha_{2}|\hat{G}_{2}^{\text{av}}(\bar{t})| \nonumber \\
&\leq -\frac{2}{\omega}\alpha_{0}(|\hat{G}_{1}^{\text{av}}(\bar{t})|+|\hat{G}_{2}^{\text{av}}(\bar{t})|)\,, \label{eq:dotVav_event_2_pf2_20240514}
\end{align}
where $\alpha_{0}=\min\{\alpha_{1}\,,\alpha_{2}\}$. From equation (\ref{eq:Vav_event_pf2}) it is easy to verify that $\sqrt{V_{\text{av}}(\bar{t})}=\sqrt{[\hat{G}_{1}^{\text{av}}(\bar{t})]^2+[\hat{G}_{2}^{\text{av}}(\bar{t})]^2}\leq \sqrt{[\hat{G}_{1}^{\text{av}}(\bar{t})]^2}+\sqrt{[\hat{G}_{2}^{\text{av}}(\bar{t})]^2}=|\hat{G}_{1}^{\text{av}}(\bar{t})|+|\hat{G}_{2}^{\text{av}}(\bar{t})|$.

Thus, $\sqrt{V_{\text{av}}(\bar{t})}\leq |\hat{G}_{1}^{\text{av}}(\bar{t})|+|\hat{G}_{2}^{\text{av}}(\bar{t})|$ and inequality (\ref{eq:dotVav_event_2_pf2_20240514}) is upper-bounded by
\begin{align}
\frac{dV_{\text{av}}(\bar{t})}{d\bar{t}}&\leq -\frac{2}{\omega}\alpha_{0}\sqrt{V_{\text{av}}(\bar{t})}\,. \label{eq:dotVav_event_2_pf2_20240514_2}
\end{align}
By defining $\tilde{V}_{\rm{av}}(\bar{t}):=\sqrt{V_{\text{av}}(\bar{t})}=\sqrt{[\hat{G}_{1}^{\text{av}}(\bar{t})]^2+[\hat{G}_{2}^{\text{av}}(\bar{t})]^2}=\|\hat{G}^{\text{av}}(\bar{t})\|$, one can write its time-derivative as $\dfrac{d\tilde{V}_{\text{av}}(\bar{t})}{d\bar{t}}=\dfrac{1}{2\sqrt{V_{\text{av}}(\bar{t})}}\dfrac{dV_{\text{av}}(\bar{t})}{d\bar{t}}$, whose upper bound is
\begin{align}
\frac{d\tilde{V}_{\text{av}}(\bar{t})}{d\bar{t}}&\leq -\frac{\alpha_{0}}{\omega}\,, \quad \bar{t}\geq 0, \label{eq:dotVav_event_2_pf2_20240514_3}
\end{align}
 by using inequality (\ref{eq:dotVav_event_2_pf2_20240514_2}). From the Comparison Theorem \cite[Theorem 7]{K:2002}, there exists an upper bound $\bar{V}_{\text{av}}(\bar{t})\geq \tilde{V}_{\text{av}}(\bar{t})$ satisfying the differential equation 
\begin{align}
\frac{d\bar{V}_{\text{av}}(\bar{t})}{d\bar{t}}=- \frac{\alpha_{0}}{\omega}\,, \quad \bar{V}_{\text{av}}(0)= \tilde{V}_{\text{av}}(0) \geq 0\,, \quad \bar{t}\geq 0. 
\end{align}
The integration of both sides of the aforementioned differential equation leads to the solution
\begin{align}
\bar{V}_{\text{av}}(\bar{t})-\bar{V}_{\text{av}}(0)=- \frac{\alpha_{0}}{\omega} \bar{t}\,, \quad \bar{t} \geq 0\,,
\end{align}
or, equivalently,
\begin{align}
\bar{V}_{\text{av}}(\bar{t})=\|\hat{G}^{\text{av}}(0)\|- \frac{\alpha_{0}}{\omega} \bar{t}\,, \quad \bar{t} \geq 0\,.
\end{align} 
Hence, there exists $\bar{t}_{1}=\frac{\omega}{\alpha_{0}}\|\hat{G}^{\text{av}}(0)\|$ such that $\bar{V}_{\text{av}}(\bar{t})= 0$, for all $t\geq t_{1}$. Now, it can be concluded that $\exists\bar{t}_{s} \in (0\,, \bar{t}_{1}]$ when the ideal sliding mode occurs, such that $\bar{V}_{\text{av}}(\bar{t})=0$, $\forall \bar{t} \geq \bar{t}_{s}$. Consequently,
\begin{align}
\|\hat{G}^{\text{av}}(\bar{t})\| &\leq \|\hat{G}^{\text{av}}(0)\|- \dfrac{\alpha_{0}}{\omega} \bar{t}\,, \quad \forall \bar{t} \in [0, \bar{t}_{s}) \,, \\ 
\|\hat{G}^{\text{av}}(\bar{t})\|&=0 \,, \quad \forall \bar{t} \geq \bar{t}_{s}\,.
\end{align}
Since $H$ is invertible and, from (\ref{eq:hatGav_event_1}), $\tilde{\theta}_{\rm{av}}(\bar{t})=H^{-1}\hat{G}_{\rm{av}}(\bar{t})$, the following inequality is established $\|H^{-1}\|^{-1}\|\tilde{\theta}_{\rm{av}}(\bar{t})\|\leq \|H\tilde{\theta}_{\rm{av}}(\bar{t})\|\leq \|H\|\|\tilde{\theta}_{\rm{av}}(\bar{t})\|$ and, therefore, we can state that
\begin{align}
\|\tilde{\theta}^{\text{av}}(\bar{t})\| &\leq \|H^{-1}\|\|H\|\|\tilde{\theta}^{\text{av}}(0)\|- \dfrac{\|H^{-1}\|\alpha_{0}}{\omega} \bar{t}\,, \quad \forall \bar{t} \in [0, \bar{t}_{s}) \,, \label{eq:tildeTheta_1_20240514} \\ 
\|\tilde{\theta}^{\text{av}}(\bar{t})\|&=0 \,, \quad \forall \bar{t} \geq \bar{t}_{s}\,. \label{eq:tildeTheta_2_20240514}
\end{align}
Since (\ref{eq:dotTildeTheta_v1_20240319}) has a discontinuous right-hand side, but is also $T$-periodic in $t$, and noting that the average system with state $\tilde{\theta}_{\text{av}}(\bar{t})$ is finite-time stable according to (\ref{eq:tildeTheta_1_20240514}) and (\ref{eq:tildeTheta_2_20240514}), we can invoke the averaging theorem in \cite[Theorem~2]{P:1979} to conclude that
\begin{align}
\|\tilde{\theta}(t)-\tilde{\theta}_{\text{av}}(t)\|\leq\mathcal{O}\left(\frac{1}{\omega}\right)\,.
\end{align}
By applying the triangle inequality \cite{A:1957}, we also obtain:
\begin{align}
&\|\tilde{\theta}(t)\|\leq\|\tilde{\theta}_{\text{av}}(t)\|+\mathcal{O}\left(\frac{1}{\omega}\right)\nonumber \\
&\leq \begin{cases} \|H^{-1}\|\|H\|\|\tilde{\theta}^{\text{av}}(0)\|\!-\! \|H^{-1}\|\alpha_{0} t\!+\!\mathcal{O}\left(\dfrac{1}{\omega}\right)\,, \forall t \in [0, t_{s})\,,  \\ \mathcal{O}\left(\dfrac{1}{\omega}\right) \,, \forall t \in [ t_{s}\,, +\infty)\,. \end{cases} \label{eq:nomrTildeTheta_pf2_gold}
\end{align}
Now, from (\ref{eq:Si}) and Fig.~\ref{fig:blockDiagram_v3}, we can verify that
\begin{align}
\theta(t)-\theta^{\ast}=\tilde{\theta}(t)+S(t)\,, \label{eq:theta_3_event_pf2}
\end{align}
where $S(t):=[S_{1}(t)\,,S_{2}(t)]^{T}$, whose the Euclidean norm satisfies
\begin{small}
\begin{align}
&\|\theta(t)-\theta^{\ast}\|=\|\tilde{\theta}(t)+S(t)\| \leq \|\tilde{\theta}(t)\|+\|S(t)\| \nonumber \\
&\mathbb{\leq} \begin{cases} \|H^{-1}\|\|H\|\|\theta(0) \mathbb{-} \theta^{\ast}\|\mathbb{-} \|H^{-1}\|\alpha_{0} t\mathbb{+}\mathcal{O}\left(a\mathbb{+}\dfrac{1}{\omega}\right)\,, \forall t \in [0, t_{s})\,,  \\ \mathcal{O}\left(a\mathbb{+}\dfrac{1}{\omega}\right) \,, \forall t \in [ t_{s}\,, +\infty)\,. \end{cases}\label{eq:theta_4_event_pf2}
\end{align}
\end{small}
leading to inequality (\ref{eq:normTheta_thm2}), for appropriate $m$ and $M_{\theta}$. \hfill $\square$

\section{SIMULATION RESULTS}

This section presents simulation results in order to illustrate the variable structure Nash equilibrium seeking based on sliding mode control. The investigated system captures a noncooperative game involving two firms operating in a duopoly market framework. These firms engage in competition aimed at maximizing their profits $J_{i}(t)$ through the pricing strategy $u_{i}(t)$ of their respective products, without sharing any information between the players.

The plant parameters are consistent with those given in \cite{FKB:2012}: initial conditions are set as $\hat{\theta}_{1}(0)=50$ and $\hat{\theta}_{2}(0)=110/3$. Additionally, $S_{d} = 100$, $p=0.20$, $m_{1}=30$ and $m_{2}=30$. According to (\ref{eq:NE_v2}), the parameter values in \cite{FKB:2012} for the duopoly game lead to a unique Nash equilibrium given by:
\begin{align}
\theta^{\ast}&=[43.3333\,, ~36.6667]^{T}\,,  \label{eq:theta_NE} \\
J^{\ast}&=[888.8889\,, ~222.2222]^{T}\,. \label{eq:J_NE}
\end{align}
The controller parameters are chosen as follows: for the NES strategy, we set $a_{1}=0.075$, $a_{2}=0.050$, $K_{1}=0.6$, $K_{2}=0.6$, $\omega_{1}=27$, and $\omega_{2}=22$. 

\begin{figure*}[h!]
	\centering
	\subfigure[Discontinuous control inputs. \label{fig:u}]{\includegraphics[width=6.5cm]{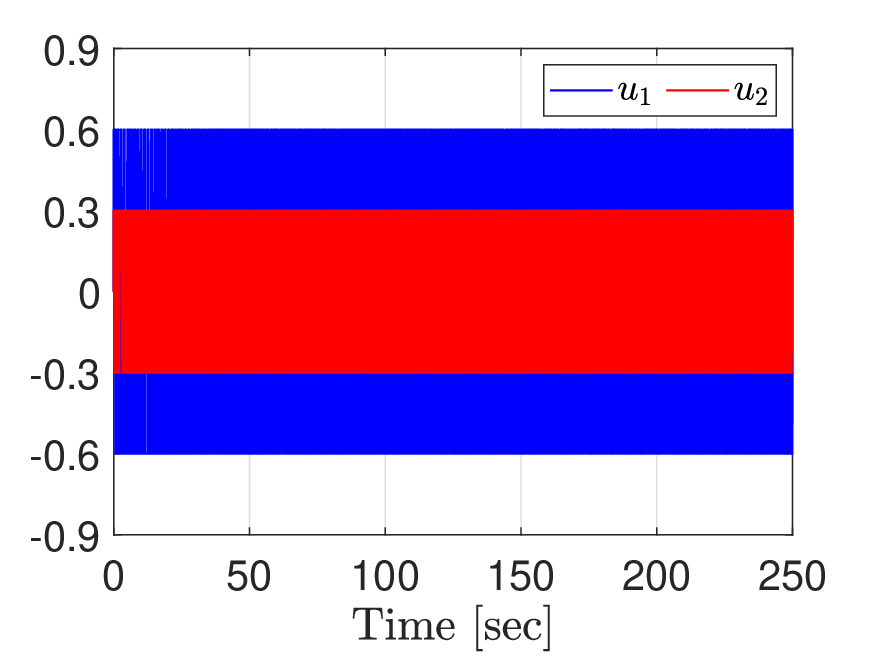}}
	~
	\subfigure[{\it Pseudogradient} estimates. \label{fig:G}]{\includegraphics[width=6.5cm]{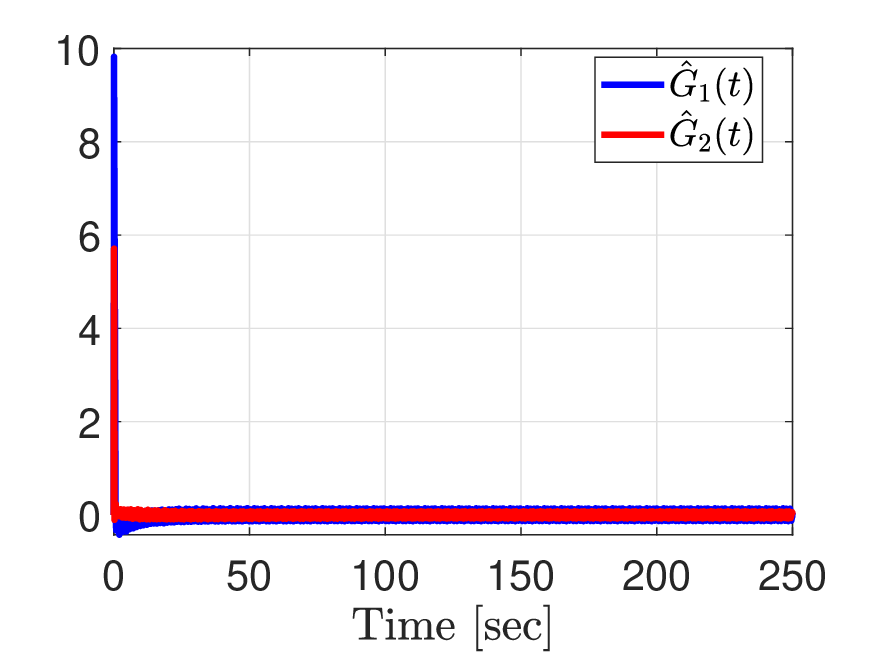}}
	\\
	\subfigure[Input signals of payoff functions, $\theta(t)$. \label{fig:theta}]{\includegraphics[width=6.5cm]{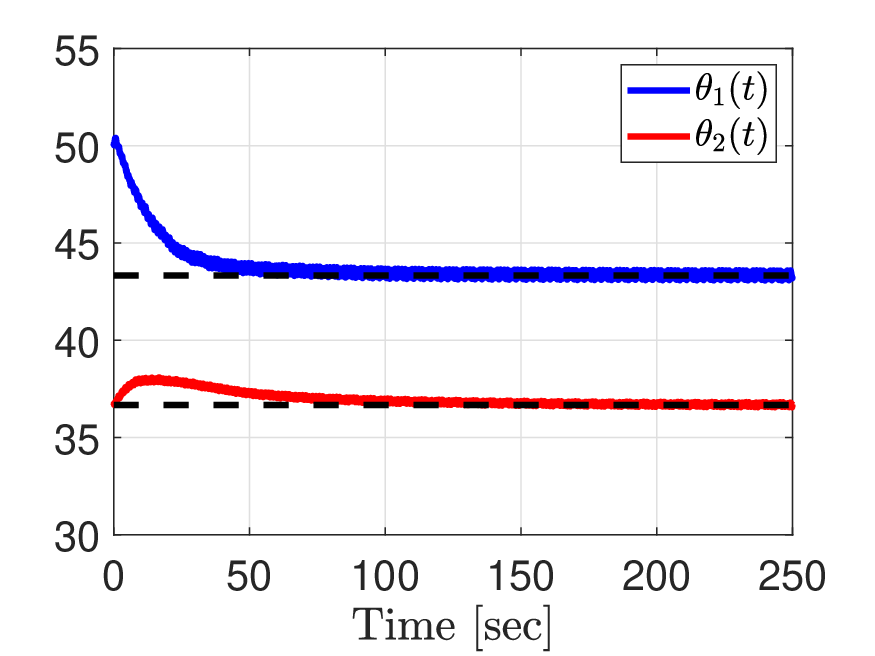}}
	~
	\subfigure[Payoff functions. \label{fig:J}]{\includegraphics[width=6.5cm]{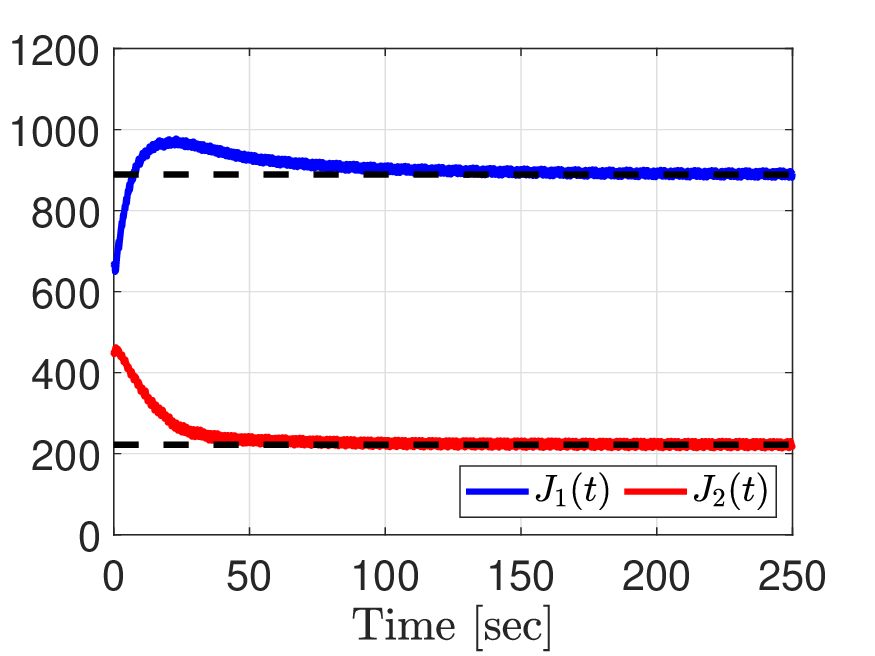}}
	\caption{Numerical simulations for NES in the noncooperative duopoly game through distributed sliding mode control policies. }\label{fig:ET_NES}
	\begin{picture}(0,0)(0,0)
		\put(35,240){\includegraphics[width=3.25cm]{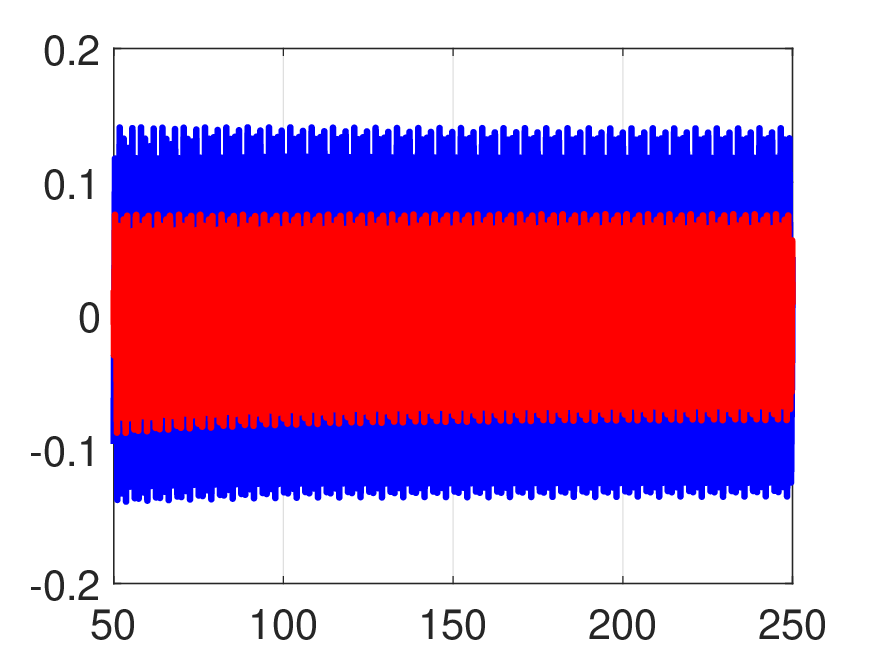}}
	\end{picture}
\end{figure*}



To achieve the Nash equilibrium (\ref{eq:theta_NE}) with the corresponding payoffs (\ref{eq:J_NE}), players P1 and P2 implement the proposed decentralized sliding mode NES strategies employing sinusoidal perturbations to execute their optimal actions.

The time evolution of the proposed approach is depicted in Fig. \ref{fig:ET_NES}. Fig. \ref{fig:u} shows the discontinuous behavior of the players' actions, consistent with the fact that each player estimates a distinct {\it pseudogradient} component and decides independently to switch its action. In a simulation spanning 150 seconds, the {\it pseudogradient} estimate is quickly brought to a neighborhood near the origin since the sliding surfaces are specified by $\hat{G}_{1}^{\text{av}}=\hat{G}_{2}^{\text{av}}=0$; see Fig.~\ref{fig:G}.

Although the  players employ distributed sliding mode strategies independently, their collective effect drives the system toward the Nash equilibrium, as depicted in Fig. \ref{fig:theta}. Additionally, Fig. \ref{fig:J} illustrates how, within the described duopoly market structure and without information sharing, each player can optimize its corresponding payoff $J_{i}(t)$ by adopting the proposed variable structure Nash equilibrium seeking via sliding mode to determine the product prices $\theta_{1}(t)$ and $\theta_{2}(t)$ in the non-cooperative game.

\section{CONCLUSIONS} 

\textcolor{black}{
This paper has presented a novel approach to attain locally stable convergence to Nash equilibrium in noncooperative games via a distributed sliding mode control scheme employing extremum seeking based on sinusoidal periodic perturbations. This methodology empowers players in a duopoly game to enhance their payoffs without necessitating information sharing, as action updates occur autonomously. Each player estimates and updates its sequence of actions based on distinct {\it pseudogradient} components by using relay functions. The proposed approach not only advances game theory but also provides a practical means for decentralized decision-making in intricate systems. Finite-time stability for the average closed-loop system is proved rather than exponential stability, generally arrived at classical extremum seeking with proportional control laws.} 

\textcolor{black}{
The results presented should be regarded as an SMC alternative to the continuous but non-smooth Nash seeking algorithms with fixed-time convergence properties in \cite{PKB:2023}.  
In the current paper, we have used discontinuous/signum/relay functions of the estimates of the {\it pseudogradient}, conducting averaging analysis of the discontinuous system with an averaging theorem of Plotnikov for differential inclusions \cite{P:1979}, and achieve finite-time convergence (with a settling time dependent on the initial condition). Specifically, the fixed-time result of \cite{PKB:2023} may look stronger since it is of fixed-time nature (the settling time is independent of the initial condition), but this conference paper provides an alternative approach, that of sliding mode, to Nash equilibrium seeking.}  


\bibliography{ifacconf}

\begin{thebibliography}{10}
\providecommand{\url}[1]{#1}
\csname url@samestyle\endcsname
\providecommand{\newblock}{\relax}
\providecommand{\bibinfo}[2]{#2}
\providecommand{\BIBentrySTDinterwordspacing}{\spaceskip=0pt\relax}
\providecommand{\BIBentryALTinterwordstretchfactor}{4}
\providecommand{\BIBentryALTinterwordspacing}{\spaceskip=\fontdimen2\font plus
\BIBentryALTinterwordstretchfactor\fontdimen3\font minus
  \fontdimen4\font\relax}
\providecommand{\BIBforeignlanguage}[2]{{%
\expandafter\ifx\csname l@#1\endcsname\relax
\typeout{** WARNING: IEEEtranS.bst: No hyphenation pattern has been}%
\typeout{** loaded for the language `#1'. Using the pattern for}%
\typeout{** the default language instead.}%
\else
\language=\csname l@#1\endcsname
\fi
#2}}
\providecommand{\BIBdecl}{\relax}
\BIBdecl

\bibitem{Sastry:2013}
S.~Amina, G.~A. Schwartz, and S.~S. Sastry, ``Security of interdependent and
  identical networked control systems,'' \emph{Automatica}, vol.~49, pp.
  186--192, 2013.

\bibitem{Aminde_SMC1}
N.~O. Aminde, T.~R. Oliveira, and L.~Hsu, ``Global output-feedback extremum
  seeking control via monitoring functions,'' in \emph{52nd IEEE Conference on
  Decision and Control (CDC)}, 2013, pp. 1031--1036.

\bibitem{Aminde_SMC2}
------, ``Multivariable extremum seeking control via cyclic search and
  monitoring function,'' \emph{International Journal of Adaptive Control and
  Signal Processing}, vol.~35, pp. 1217--1232, 2021.

\bibitem{A:1957}
T.~Apostol, \emph{Mathematical Analysis - A Modern Approach to Advanced
  Calculus}.\hskip 1em plus 0.5em minus 0.4em\relax Addison-Wesley Publishing
  Company, 1957.

\bibitem{B:1987}
T.~Ba\c{s}ar, ``Relaxation techniques and asynchronous algorithms for on-line
  computation of noncooperative equilibria,'' \emph{J. Economic Dynamics and
  Control}, vol.~11, pp. 531--549, 1987.

\bibitem{Basar:1999}
T.~Ba\c{s}ar and G.~J. Olsder, \emph{Dynamic Noncooperative Game Theory}.\hskip
  1em plus 0.5em minus 0.4em\relax SIAM Series in Classics in Applied
  Mathematics, 1999.

\bibitem{BZ:2018}
T.~Ba\c{s}ar and G.~Zaccour, Eds., \emph{Handbook of Dynamic Game Theory,
  Volume I (Theory of Dynamic Games)}.\hskip 1em plus 0.5em minus 0.4em\relax
  Springer International Publishing, 2018.

\bibitem{BZ2:2018}
------, \emph{Handbook of Dynamic Game Theory, Volume II (Applications of
  Dynamic Games)}.\hskip 1em plus 0.5em minus 0.4em\relax Springer
  International Publishing, 2018.

\bibitem{FKB:2012}
P.~Frihauf, M.~Krsti{\' c}, and T.~Ba{\c s}ar, ``Nash equilibrium seeking in
  noncooperative games,'' \emph{IEEE Trans. Automat. Contr.}, vol.~57, pp.
  1192--1207, 2012.

\bibitem{Tirole:1991}
D.~Fudenberg and J.~Tirole, \emph{Game Theory}.\hskip 1em plus 0.5em minus
  0.4em\relax The MIT Press, 1991.

\bibitem{GKN:2012}
A.~Ghaffari, M.~Krsti{\' c}, and D.~Ne{\u s}ic, ``Multivariable {N}ewton-based
  extremum seeking,'' \emph{Automatica}, vol.~48, pp. 1759--1767, 2012.

\bibitem{Basar:2019}
Z.~Han, D.~Niyato, W.~Saad, and T.~Ba{\c s}ar, \emph{Game Theory for Next
  Generation Wireless and Communication Networks: Modeling, Analysis, and
  Design}.\hskip 1em plus 0.5em minus 0.4em\relax Cambridge University Press,
  2019.

\bibitem{horn1985}
R.~A. Horn and C.~R. Johnson, \emph{Matrix Analysis}.\hskip 1em plus 0.5em
  minus 0.4em\relax Cambridge, U.K.: Cambridge Univ. Press, 1985.

\bibitem{K:2002}
H.~K. Khalil, \emph{Nonlinear Systems}.\hskip 1em plus 0.5em minus 0.4em\relax
  Prentice Hall, 2002.

\bibitem{K:2014}
M.~Krsti{\' c}, ``Extremum seeking control,'' in \emph{Encyclopedia of Systems
  and Control}.\hskip 1em plus 0.5em minus 0.4em\relax London: Springer, 2014,
  vol.~1, pp. 413--416.

\bibitem{KW:2000}
M.~Krsti{\' c} and H.-H. Wang, ``Stability of extremum seeking feedback for
  general nonlinear dynamic systems,'' \emph{Automatica}, vol.~36, pp.
  595--601, 2000.

\bibitem{LB:1987}
S.~Li and T.~Ba\c{s}ar, ``Distributed learning algorithms for the computation
  of noncooperative equilibria,'' \emph{Automatica}, vol.~23, pp. 523--533,
  1987.

\bibitem{M:2011}
\BIBentryALTinterwordspacing
J.~A. Moreno, ``Lyapunov approach for analysis and design of second order
  sliding mode algorithms,'' in \emph{Sliding Modes after the First Decade of
  the 21st Century}, ser. Lecture Notes in Control and Information Sciences,
  L.~Fridman, J.~Moreno, and R.~Iriarte, Eds.\hskip 1em plus 0.5em minus
  0.4em\relax Berlin, Heidelberg: Springer, 2011, vol. 412, pp. 113--149.
  [Online]. Available: \url{https://doi.org/10.1007/978-3-642-22164-4_4}
\BIBentrySTDinterwordspacing

\bibitem{MO:2012}
J.~A. Moreno and M.~Osorio, ``Strict lyapunov functions for the super-twisting
  algorithm,'' \emph{IEEE Transactions on Automatic Control}, vol.~57, no.~4,
  pp. 1035--1040, 2012.

\bibitem{N:1951}
J.~F. Nash, ``Noncooperative games,'' \emph{Annals of Mathematics}, vol.~54,
  pp. 286--295, 1951.

\bibitem{Oliveira_SMC1}
T.~R. Oliveira, L.~Hsu, and A.~J. Peixoto, ``Output-feedback global tracking
  for unknown control direction plants with application to extremum-seeking
  control,'' \emph{Automatica}, vol.~47, pp. 2029--2038, 2011.

\bibitem{Oliveira_SMC2}
T.~R. Oliveira, A.~J. Peixoto, and L.~Hsu, ``Global real-time optimization by
  output-feedback extremum-seeking control with sliding modes,'' \emph{Journal
  of the Franklin Institute}, vol. 349, pp. 1397--1415, 2012.

\bibitem{Ozguner_SMC}
Y.~Pan, {\"U}.~{\"O}zg{\"u}ner, and T.~Acarman, ``Stability and performance
  improvement of extremum seeking control with sliding mode,''
  \emph{International Journal of Control}, vol.~76, pp. 968--985, 2003.

\bibitem{P:1979}
V.~A. Plotnikov, ``Averaging of differential inclusions,'' \emph{Ukrainian
  Mathematical Journal}, vol.~31, pp. 454--457, 1980.

\bibitem{PKB:2023}
J.~I. Poveda, M.~Krsti{\' c}, and T.~Ba\c{s}ar, ``Fixed-time {Nash} equilibrium
  seeking in time-varying networks,'' \emph{IEEE Trans. Automat. Contr.},
  vol.~68, pp. 1954--1969, 2023.

\bibitem{R:1965}
J.~B. Rosen, ``Existence and {U}niqueness of {E}quilibrium {P}oints for
  {C}oncave {N}-{P}erson {G}ames,'' \emph{Econometrica}, vol.~33, no.~3, pp.
  520--534, 1965.

\bibitem{taussky1949}
O.~Taussky, ``A recurring theorem on determinants,'' \emph{Amer. Math.
  Monthly}, vol.~56, no.~10, pp. 672--676, 1949.

\bibitem{TK:2023}
V.~Todorovski and M.~Krsti{\' c}, ``Practical prescribed-time seeking of a
  repulsive source by unicycle angular velocity tuning,'' \emph{Automatica},
  vol.~68, p. 111069, 2023.

\bibitem{YK:2022}
C.~T. Yilmaz and M.~Krsti{\' c}, ``Prescribed-time extremum seeking with chirpy
  probing for {PDEs--Part I: Delay},'' \emph{IEEE American Control Conference},
  pp. 1000--1005, 2022.

\bibitem{ZTB:2013}
Q.~Zhu, H.~Tembine, and T.~Ba\c{s}ar, ``Hybrid learning in stochastic games and
  its applications in network security,'' in \emph{Series on Computational
  Intelligence}, F.~L. Lewis and D.~Liu, Eds.\hskip 1em plus 0.5em minus
  0.4em\relax IEEE Press/Wiley, 2013, ch.~14, pp. 305--329.

\end{thebibliography}

\end{document}